\input amstex
\documentstyle{amsppt}
\magnification=\magstep1 \NoRunningHeads
\NoBlackBoxes

\topmatter
\title
Flows with uncountable but meager \\ group of self-similarities
\endtitle
\author
Alexandre I. Danilenko
\endauthor
\address
 Institute for Low Temperature Physics
\& Engineering of National Academy of Sciences of Ukraine, 47 Lenin Ave.,
 Kharkov, 61164, UKRAINE
\endaddress
\email alexandre.danilenko\@gmail.com
\endemail

\dedicatory
Dedicated to A. M. Stepin on the occasion of his 70-th birthday
\enddedicatory

\abstract
Given an ergodic probability preserving flow $T=(T_t)_{t\in\Bbb R}$,
let $I(T):=\{s\in\Bbb R^*\mid T\text{ is isomorphic to }(T_{st})_{t\in\Bbb R}\}$.
A weakly mixing Gaussian flow $T$ is constructed such that $I(T)$ is uncountable and meager.
For a Poisson flow $T$, a subgroup $I_{\text{Po}}(T)\subset I(T)$ of Poissonian self-similarities is introduced.
Given a probability measure $\kappa$ on $\Bbb R^*_+$, a  Poisson flow $T$ is constructed such that $I_{\text{Po}}(T)$ is the  group of $\kappa$-quasi-invariance.
\endabstract
\endtopmatter

\document

\head 0. Introduction
\endhead

Throughout the paper $\Bbb R^*$ denotes the multiplicative group of reals different from $0$ and  $\Bbb R^*_+$ denotes the multiplicative group of positive reals.
Let $T=(T_t)_{t\in\Bbb R}$ be an ergodic free measure preserving flow on a standard non-atomic probability space $(X,\goth B,\mu)$.
Given $s\in\Bbb R^*$, we denote by $T\circ s$ the  flow $(T_{st})_{t\in\Bbb R}$.
Let
$$
I(T):=\{s\in\Bbb R^*\mid T\circ s\text{ is isomorphic to }T\}.
$$
It is easy to see that  $I(T)$ is a multiplicative subgroup of $\Bbb R^*$.
It is called the group of {\it self-similarities} of $T$.
If $I(T)\ne\{1,-1\}$  then $T$ is called {\it  self-similar}.
There are a lot of problems related to the self-similarities of ergodic flows.
 We refer the reader to  the recent papers \cite{FrLe}, \cite{DaRy}  and references therein for details.
In particular, it was shown there that given any countable subgroup $G$ in $\Bbb R$, there is a weakly mixing flow $T$ with $I(T)=G$.
In these notes we consider Problem 2 from the ``Open problems'' section in \cite{FrLe}:
\roster
\item"---"
is there an ergodic flow $T$ for which the group $I(T)$ is uncountable but has zero Lebesgue measure?
\endroster
Examples of such flows were constructed recently in \cite{DaRy, Section~2}.
However the construction there incorporates essentially some subtle facts from the measurable orbit theory such as the theorem on outer conjugacy for groups of automorphisms of continuous ergodic equivalence relations \cite{VeFe} and Ratner's theory on joinings of horocycle flows \cite{Ra}.
Our purpose here is to provide  other (simpler) examples, independent of \cite{VeFe} and \cite{Ra}.

\proclaim{Theorem 0.1}
There exist weakly mixing Gaussian flows with uncountable but meager group of self-similarities.
\endproclaim

In this connection we recall that
\roster
\item"---"
$I(T)$ is a Borel subset of $\Bbb R^*$ for each flow $T$ \cite{DaRy} and
\item"---"
a Borel subgroup of $\Bbb R^*$ is meager if and only if it has zero Lebesgue measure.
\endroster
It should be noted that the flows constructed in Theorem~0.1 are different (in fact, disjoint in the sense of Furstenberg \cite{Fu}) from the mixing flows with uncountable meager group of self-similarities constructed in \cite{DaRy, Section 2}.

Next, for each Poisson flow $T$, we define a subgroup  $I_{\text{Po}}(T)\subset I(T)$
of Poissonian self-similarities.
It consists of those $s\in\Bbb R^*$ such that $T\circ s$ is conjugate to $T$ via a Poisson transformation.
Given a probability measure $\tau$ on a locally compact second countable  Abelian group $G$,
 let $H_G(\tau)$ be the set of all $g\in G$ such the translation of $\tau$ by $g$ is equivalent to $\tau$.
Then $H_G(\tau)$ is a Borel subgroup in $G$ (see \cite{Na} and references therein).
It is called the {\it group of $\tau$-quasi-invariance}.

Recently, Ryzhikov  \cite{Ry} constructed  a class of  Poisson flows $T$ with {\it total self-similarity}, i.e. $I(T)=\Bbb R^*$\footnote{The Bernoulli flow with infinite entropy and horocycle flows are other examples of flows with total self-similarities.}.
Generalizing his construction,
we prove the following claim.

\proclaim{Theorem 0.2}
For each probability measure $\kappa$ on $\Bbb R^*_+$, there is a  weakly mixing Poisson flow $\widetilde T$  with $I_{\text{Po}}(\widetilde T) =H_{\Bbb R^*_+}(\kappa)$.
\endproclaim

The author thanks E. Roy and the anonymous referee for the useful remarks.

\head  1. Gaussian examples
\endhead

 Fix an onto homomorphism  $\pi:\Bbb R\to\Bbb T$.
Then  for each
 probability measure $\rho$ on $\Bbb T$, there is
a probability measure $\rho'$ on $\Bbb R$ such that $\rho'\circ\pi^{-1}\sim\rho$ and
$H_\Bbb R(\rho'):=\pi^{-1}(H_\Bbb T(\rho))$.
We call $\rho'$ a {\it standard lift} of $\rho$ to $\Bbb R$.
It is easy to verify that the equivalence class of the standard lift of $\rho$ is defined
uniquely by the equivalence class of $\rho$.

Consider an isomorphism $\vartheta:\Bbb R\to\Bbb R_+^*$.
Let $\sigma$ denote the only non-atomic symmetric probability measure on $\Bbb R$
such that $\sigma\restriction\Bbb R_+^*=\frac 12\rho'\circ\vartheta^{-1}$.
We use here a natural (set theoretical only, not group theoretical) embedding $\Bbb R_+^*\subset\Bbb R$.
Then
$$
H_{\Bbb R^*}(\sigma)=\vartheta(H_{\Bbb R}(\rho'))\cup(- \vartheta(H_{\Bbb R}(\rho'))).
$$
It follows that if $H_{\Bbb T}(\rho)$ is uncountable and meager in $\Bbb T$ then $H_{\Bbb R^*}(\sigma)$ is uncountable and meager in $\Bbb R^*$.

Suppose that $\rho$ is non-atomic and ergodic with respect to the non-singular action of $H_{\Bbb T}(\rho)$ on $(\Bbb T,\rho)$ by translations.
We recall (see \cite{AaNa} and \cite{Na}) that $H_{\Bbb T}(\rho)$ admits a unique natural topology in which it is  a Polish group.
Then $\rho$ is called $H_{\Bbb T}(\rho)$-ergodic if it is $G$-ergodic
for some (and hence for every) dense countable subgroup $G$ of $H_{\Bbb T}(\rho)$.
It follows that
 $\rho'$ is  $H_{\Bbb R}(\rho')$-ergodic.
We also note that $H_{\Bbb R}(\rho')$ is dense in $\Bbb R$.
Next,  $\sigma$ is ergodic with respect
$H_{\Bbb R^*}(\sigma)$
acting on $\Bbb R^*$ by multiplication
and $H_{\Bbb R^*}(\sigma)$
is dense in $\Bbb R^*$.

Let $T$ denote the Gaussian flow associated with $\sigma$ (see \cite{Co-Si}).

\proclaim{Proposition 1.1}
If $\rho$ is a non-atomic probability measure on $\Bbb T$ such that
\roster
\item
$H_{\Bbb T}(\rho)$ is uncountable and meager and
\item
$\rho$ is ergodic with respect to
$H_{\Bbb T}(\rho)$
\endroster
then  $I(T)\supset H_{\Bbb R^*}(\sigma)$  but $\Bbb R^*_+\not\subset
I(T)$.
\endproclaim

\demo{Proof}
Since $T$ is a Gaussian flow associated with $\sigma$, it follows that $I(T)\supset H_{\Bbb R^*}(\sigma)$ \cite{FrLe}.
On the other hand, it is well known that
the
  reduced maximal spectral type of $T$
is given by
the measure  $\exp'\sigma:=\sum_{p\ge 1}\frac{\sigma^{*p}}{p!}$, where $\sigma^{*p}$  denotes the $p$-th convolution
power of $\sigma$ \cite{Co-Si}.
Hence $I(T)\subset H_{\Bbb R^*}(\exp'\sigma)$.
We claim that
 $I(T)\ne \Bbb R^*$.
Indeed, otherwise
$H_{\Bbb R^*}(\exp'\sigma)=\Bbb R^*$ and hence
$\exp'\sigma$ is equivalent to Lebesgue measure $\lambda_{\Bbb R^*}$.
Since $\lambda_{\Bbb R^*}$ and $\sigma$ are both ergodic with respect to $H_{\Bbb R^*}(\sigma)$ and $\lambda_{\Bbb R^*}\sim\exp'\sigma\succ\sigma$, it follows that $\lambda_{\Bbb R^*}\sim\sigma$
and hence $H_{\Bbb R^*}(\sigma)=\Bbb R^*$, a contradiction.
In a similar way we can verify that
$I(T)\not\supset\Bbb R^*_+$.
\qed
\enddemo

\demo{Proof of Theorem~0.1}
Examples of $\rho$ satisfying (1) and (2) are given in \cite{AaNa, Theorem~4.1} and  \cite{Pa}.
Now it suffices to note that
\roster
\item"---" $T$ is weakly mixing whenever $\sigma$ is non-atomic,
 \item"---"
if a subgroup of $\Bbb R^*$ does not contain $\Bbb R^*_+$ then it is meager in $\Bbb R^*$
\endroster
and apply Proposition~1.1.
\qed
\enddemo

We now briefly explain how to construct measures satisfying (1) and (2).
Suppose we have a sequence of positive integers $(n_j)_{j\ge 1}$
and a sequence of complex numbers $(a_j)_{j\ge 1}$
such that $n_j>2(n_1+\dots+n_{j-1})$  and $|a_j|\le 1$ for all $j$.
We let $P_j(z):=1+a_jz^{n_j}/2+\overline{a_j}z^{-n_j}/2$ for $z\in\Bbb T$.
Then the polynomial $P_j$ is non-negative and $\int_{\Bbb T}\prod_{k=1}^j P_k\,d\lambda_\Bbb T=1$ for each $j$, where $\lambda_\Bbb T$ is the normed Haar measure on $\Bbb T$.
The sequence of probability measures $(\prod_{k=1}^jP_k)\lambda_{\Bbb T}$
converges weakly as $j\to\infty$ to a probability measure $\rho$ on $\Bbb T$ which is called the {\it Riesz product} associated with $(n_j)_{j\ge 1}$ and $(a_j)_{j\ge 1}$ \cite{Na}.
The Fourier transform of $\rho$ is explicitly known: $\widehat\rho(n)=0$
except if $n$ is a finite sum
$\sum_{j}k_jn_j$ with $k_j\in\Bbb Z$ and $|k_j|\le 1$\footnote{Due to the condition imposed on $(n_j)_{j\ge 1}$, if $n$ can be written as $\sum_{j}k_jn_j$ then this representation is unique. Hence $\widehat\rho(n)$ is well defined.},
 in which case
$\widehat\rho(n)=\prod_j a_j^{k_j}$.
Using standard methods of harmonic analysis and classical facts about
equivalence of Riesz products,
F.~Parreau shows in \cite{Pa} that
\roster
\item"---"
$
H(\rho)=\{z\in\Bbb T\mid \sum_{j\ge 1}|a_j|^2|1-a_jz^{n_j}|^2<\infty\}
$
and
\item"---"
 if the lacunary condition
$
\sum_{j\ge 1}|a_j|^2(n_j/n_{j+1})^2<\infty
$
 holds then
$\rho$ is ergodic with respect to $H(\rho)$.
Moreover, in this case $H(\rho)$ is uncountable.
\item"---" If $\sum_i|a_j|^2=\infty$ then $H(\rho)\ne\Bbb T$.
\endroster
Thus we see that if $\rho$ is the Riesz product associated with $n_j=j!$ and $a_j= 1$ for all $j$ then $\rho$ satisfies (1) and (2).

\remark{Remark 1.2} We note that the flows constructed in Proposition~1.1 are quite different from flows with uncountable but meager group of self-similarities constructed in \cite{DaRy, \S3}.
To explain this precisely we recall some concepts from the theory of joinings of dynamical systems (see \cite{Th}, \cite{dR}).
Given two ergodic flows $(X,\mu, T)$ and $(Y,\nu,S)$, a joining of $T$ and $S$ is a $T\times S$-invariant measure on $X\times Y$ whose marginals on $X$ and $Y$ are $\mu$ and $\nu$ respectively.
For instance, $\mu\times\nu$ is a joining of $T$ and $S$.
If there are no other joinings then $T$ and $S$ are called disjoint \cite{Fu}.
Of course, disjoint flows are non-isomorphic.
A flow $(X,\mu,T)$ is called 2-fold quasi-simple if for every ergodic joining $\lambda\ne\mu\times\mu$ of $T$ with itself, the two marginal projections $(X\times X,\lambda)\to(X,\mu)$
are finite-to-one $\mod \lambda$.\footnote{More generally, the flow is called 2-fold quasi-simple  if the aforementioned marginal extensions are isometric \cite{RyTh} or even distal \cite{dJLe}.}
The examples from \cite{DaRy, \S 3} are 2-point extensions of some horocycle flows.
Hence they are quasi-simple according to \cite{RyTh} and \cite{dJLe}.
By \cite{RyTh} and \cite{dJLe}, a weakly mixing Gaussian flow is disjoint with  each quasi-simple flow.
Thus the examples constructed in Theorem~0.1 are disjoint with the flows constructed in \cite{DaRy, \S 3}.
\endremark

\head
2. Poisson flows
\endhead

Let $X$ be a locally compact non-compact second countable space.
Denote by $\widetilde X$ the vector space of Radon measures on $X$.
We endow $\widetilde X$ with the $*$-weak topology, i.e.
 the weak topology generated by the natural duality of $\widetilde X$ and the
space $C_0(X)$ of continuous functions with compact support on $X$.
The weak topology on $\widetilde X$ generated by this duality is called the
 $*$-weak topology.
Since the Borel structure $\widetilde{\goth B}$ generated by the $*$-weak topology is the same as the Borel structure generated by the strong topology on $\widetilde X$ and the strong topology is Polish, it follows that
$\widetilde {\goth B}$ is standard Borel.
We fix a homeomorphism $T$ of $X$.
It induces canonically a homeomorphism $\widetilde T$ of $\widetilde X$.
Denote by $\Cal M_{+,T}(X)\subset\widetilde X$ the closed cone of non-negative ($\sigma$-finite) Radon $T$-invariant measures  and denote by $\Cal M_{+,1,\widetilde T}(\widetilde X)$ the set of probability $\widetilde T$-invariant measures on $(\widetilde X,\widetilde {\goth B})$.
Given a non-negative measure $\mu$ on $X$,
a probability measure $\widetilde\mu$ on $\widetilde X$ is well defined by the following two conditions:
\roster
\item"(c1)"
$\widetilde\mu([K,j])=\exp(-\mu(K))\mu(K)^j/j!$ for each  compact subset $K\subset X$ and each non-negative integer $j$,
\item"(c2)"
$\widetilde\mu([K,j]\cap[K',j'])=
\widetilde\mu([K,j])\widetilde\mu([K',j'])$ whenever $K\cap K'=\emptyset$,
 \endroster
where $[K,j]$ stands for the cylinder $\{\omega\in\widetilde X\mid
\omega(K)=j\}$.
If $\mu$ is $T$-invariant then $\widetilde \mu$ is $\widetilde T$-invariant.
The dynamical system $(\widetilde X,\widetilde\mu,\widetilde T)$ is called the {\it Poisson suspension} of $(X,\mu,T)$.
This definition can be naturally extended to the case where $(X,\mu, T)$
is an arbitrary  standard Borel $\mu$-preserving dynamical system with $\mu$ infinite and $\sigma$-finite (see, e.g., \cite{Ro}, \cite{Ja--dR}).
If a probability preserving transformation is measure theoretically isomorphic to a Poisson suspension then it is called a {\it Poisson} transformation.

We recall an interesting example by Ryzhikov  from his recent work \cite{Ry}.
 Let $V=(V_t)_{t\in\Bbb R}$ be an ergodic  conservative measure preserving flow on an infinite measure space $(Z,\nu)$.
We now define a conservative measure preserving flow $T=(T_t)_{t\in\Bbb R}$ on the infinite measure space $(X,\mu):=(\Bbb R^*\times Z,\lambda_{\Bbb R^*}\times\nu)$ by
$$
T_t(s,z):=(s, V_{st}z).
$$
Since $T$ has no invariant subsets of finite measure, $\widetilde T$ is weakly mixing \cite{Ro}.
It is easy to see that $I(T)=\Bbb R^*$.
 The self-similarity group for the infinite measure preserving flows is defined in the very same way as for the probability preserving flows.
Hence for the Poisson suspension $\widetilde T=(\widetilde {T_t})_{t\in\Bbb R}$ of $T$, we also obtain
 $I(\widetilde T)=\Bbb R^*$.

We now consider a generalization of Ryzhikov's example (see the construction in the proof of Theorem~0.2).
Let $\widetilde T=(\widetilde T_t)_{t\in\Bbb R}$ be a  Poisson suspension of a measure preserving flow  $T=(T_t)_{t\in\Bbb R}$ on a
$\sigma$-finite measure space $(X,\mu)$.
We let
$
I_{\text{Po}}(\widetilde T):=I(T)$.
Then $I_{\text{Po}}(\widetilde T)$ is obviously a subgroup of $I(\widetilde T)$.
We call it the {\it group of Poissonian self-similarities} of $\widetilde T$.

\demo{Proof of Theorem 0.2} Let $V$ be an ergodic  conservative measure preserving flow $V$ on an infinite  $\sigma$-finite measure space $(Y,\nu)$.
 Let $\sigma$ stand for a measure of the maximal spectral type of $V$.
For $s>0$, we denote by $\sigma_s$ the image of $\sigma$ under
the mapping $\Bbb R\ni t\mapsto ts\in\Bbb R$.
Then $\sigma_s$ is a measure of the maximal spectral type of $V\circ s$.
Suppose that  $\sigma_s\perp\sigma$ for each $1\ne s>0$ (see \cite{DaRy} for examples of such flows).
Fix a probability  measure $\kappa$ on $\Bbb R^*_+$.
We now define a measure preserving flow $T=(T_t)_{t\in\Bbb R}$ on the infinite measure space $(X,\mu):=(\Bbb R_+^*\times\Bbb R\times Z,\kappa\times\lambda_{\Bbb R}\times\nu)$ by
$$
T_t(s,y,z):=(s,y, V_{st}z).
$$
It follows from our assumption on $V$ that if $V\circ s=RV\circ {s'}R^{-1}$
for  $s,s'>0$ and a measure preserving transformation $R$ of $(Y,\nu)$ then $s=s'$.

Since $T$ has no invariant subsets of finite measure,
it follows that the Poisson suspension $\widetilde T$ of $T$
is weakly mixing.
We claim that $I_{\text{\rom{Po}}}(\widetilde T)=H_{\Bbb R^*_+}(\kappa)$.
It is easy to see that $H_{\Bbb R^*_+}(\kappa)\subset I(T)$.
Indeed, if $h\in H_{\Bbb R^*_+}(\kappa)$ then we let
$$
Q(s,y,z):=(hs,y\cdot\frac{d\kappa}{d\kappa_h}(s),z).
$$
Then $Q$ is a $\mu$-preserving  transformation of $X$ and
$Q^{-1}TQ=T\circ h$.
Now we prove the converse inclusion $I(T)\subset H_{\Bbb R^*_+}(\kappa)$.
Take $h\in I(T)$.
Then there  is a measure preserving transformation $Q$ of $X$
such that $Q^{-1}TQ=T\circ h$.
It is easy to write the $T$-ergodic decomposition of $\mu$:
$$
\mu=\int_{\Bbb R^*_+\times \Bbb R}\delta_s\times\delta_y\times\nu \, d\kappa(s)dy
\tag2-1
$$
The restriction of $T$ to $(X,\delta_s\times\delta_y\times\nu)$ is isomorphic to $V\circ s$.
Therefore two ergodic components
$(X,\delta_s\times\delta_y\times\nu,T)$ and $(X,\delta_{s'}\times\delta_{y'}\times\nu,T)$
are isomorphic if and only if $s=s'$.
We note that \thetag{2-1} is also the  $T\circ h$-ergodic decomposition
of $\mu$.
However the ergodic component
$(X,\delta_s\times\delta_y\times\nu,T\circ h)$ is isomorphic to $V\circ (hs)$.
Since $Q$ conjugates each ergodic component of $T$ with an ergodic component of $T\circ h$, we obtain
that
$$
Q(s,y,z)=(h^{-1}s,\bullet, \bullet).
$$
Thus the marginal projection $X\ni(s,y,z)\mapsto s\in\Bbb R^*_+$ intertwines $Q$ with the mapping  $\Bbb R^*_+\ni s\mapsto h^{-1}s\in\Bbb R^*_+$.
Since $Q$ is $\mu$-preserving, it follows that this mapping is $\kappa$-nonsingular, i.e.
 $h\in H_{\Bbb R^*_+}(\kappa)$.
\qed
\enddemo

\head 3. Concluding remarks and problems
\endhead

\roster
\item (Ryzhikov's question)
Are there non-mixing weakly mixing   flows with uncountable group of self-similarities?
The examples from \cite{DaRy, Section 2} are mixing.
We conjecture that the class of Gaussian flows constructed in Theorem~0.1 contains non-mixing flows.
\item Are there {\it prime} weakly mixing flows with uncountable but meager group of of self-similarities?
 The flows constructed
 in  \cite{DaRy, Section~2} and Theorem~0.1 have non-trivial factors.
\item
We do not know whether  $I_{\text{Po}}(\widetilde T)=I(\widetilde T)$ for the Poisson flows $\widetilde T$ constructed in Theorem~0.2.
Are there some extra conditions on $\kappa$ and $V$ which imply the equality $I_{\text{Po}}(\widetilde T)=I(\widetilde T)$?
It is easy to deduce from \cite{Ro, Proposition~5.2} that this equality holds if
the maximal spectral type $\tau$ of $T$ is orthogonal to $\sum_{j>1}\frac 1{j!}\tau^{*j}$.
A direct calculation shows that $\tau=\int_{\Bbb R^*}\sigma_s\,d\kappa(s)$, where $\sigma$ is a measure of the maximal spectral type of $V$.
\item
On the other hand, it is also interesting to find a Poisson flow $\widetilde T$ with
$I_{\text{Po}}(\widetilde T)\ne I(\widetilde T)$.

\endroster

\enddocument